
\magnification=\magstep1
\hsize=15truecm
\tolerance=4000
\vsize=22truecm
\mathsurround=1pt
\nopagenumbers

\voffset=2\baselineskip
\hoffset=.3truecm

\def\highskip{\medskip}
\def\highlight{\leftskip=50pt
               \parindent=0pt}
\def\endhighlight{
                  \leftskip=0pt
                  \parindent=20pt}

\def\abstract{\leftskip=20pt
              \rightskip=20pt
              \parindent=0pt}

\def\endabstract{\leftskip=0pt
              \rightskip=0pt
              \parindent=20pt}

\font\bigbf=cmbx10 scaled \magstep1

\font\bigbigbigbf=cmbx10 scaled \magstep3

\font\bigrm=cmr10 scaled \magstep1

\def\chapter#1{\par\bigbreak\smallskip \centerline{\bigbf #1}\medskip}

\def\section#1{\par\bigbreak {\bf #1}\nobreak\enspace}

\def\k{\kappa}
\def\e{\epsilon}
\def\o{\omega}
\def\s{\sigma}
\def\d{\delta}
\def\f{\varphi}
\def\l{\lambda}
\def\r{\rho}
\def\a{\alpha}
\def\b{\beta}
\def\p{\pi}
\def\n{\eta}
\def\z{\zeta}
\def\t{\tau}
\def\g{\gamma}
\def\h{\theta}
\def\x{\xi}

\def\G{\Gamma}

\def\m{\mu}




\def\lause #1. #2\par\par{\medbreak\smallskip
\smallskip\smallskip
\noindent{\bf#1.
\enspace}{\sl#2\par}\par\medbreak}

\def\sublemma #1. #2\par\par{\medbreak\smallskip
\noindent{\bf#1.
\enspace}{\sl#2\par}\par\medbreak}

\def\maar #1 #2 {\par\par
\smallskip\smallskip
\smallskip\medbreak\noindent{\bf#1 #2
\enspace}}

\def\tod #1 {\par\par\medbreak\noindent{\bf #1 \enspace}}

\def\ceisi #1 {\par\par\medbreak\noindent{\bf #1}\ \ }

\def\valit {\par\par
\smallskip\smallskip
\smallskip\medbreak}

\def\rem #1 #2. #3\par{\medbreak{\bf #1 #2.
\enspace}{#3}\par\medbreak}

\def\sqr#1#2{{\vcenter{\hrule height.#2pt
      \hbox{\vrule width.#2pt height#1pt \kern#1pt
         \vrule width.#2pt}
       \hrule height.#2pt}}}
\def\eop{\mathchoice\sqr34\sqr34\sqr{2.1}3\sqr{1.5}3}

                                                                     %
                                                                     %
\newdimen\refindent\newdimen\plusindent                              %
\newdimen\refskip\newdimen\tempindent                                %
\newdimen\extraindent                                                %
                                                                     %
                                                                     %
\def\ref#1 #2\par{\setbox0=\hbox{#1}\refindent=\wd0                  %
\plusindent=\refskip                                                 %
\extraindent=\refskip                                                %
\advance\extraindent by 30pt                                         %
\advance\plusindent by -\refindent\tempindent=\parindent 
\parindent=0pt\par\hangindent\extraindent 
{#1\hskip\plusindent #2}\parindent=\tempindent}                      %
\refskip=1.5cm                                               %
                                                                     %


\def\A{{\cal A}}
\def\B{{\cal B}}

\def\M{{\cal M}}

\def\k{\kappa}

\def\yesheadline{\headline={\hfil\tenrm\folio}}

\yesheadline

\def\ti{\item{(i)}}
\def\tii{\item{(ii)}}
\def\tiii{\item{(iii)}}
\def\tiv{\item{(iv)}}

\def\tar{\models}

\def\kuva#1 #2{
    \midinsert
    \vskip #2
    \centerline{#1}
    \endinsert}

\def\eks{\exists}
\def\kai{\forall}

\def\o{\omega}

\def\eset{\emptyset}

\def\itskip{\smallskip}

\def\loppu{\vfill\eject\end}

\def\lo1o{$L_{\o_1\o}$}

\def\inv{^{-1}}

\def\Fraisse{Fra\"\i ss\'e }

\def\dv{\vert \vert}
\def\ve{\vert}
\def\raj{\mid}

\def\dom{{\rm dom}}
\def\ran{{\rm ran}}

\def\max{{\rm max}}

\def\height{{\rm height}}
\def\AUT{{\rm AUT}}

\yesheadline

\bigskip
\centerline{\bigbigbigbf On the number of automorphisms}
\centerline{\bigbigbigbf of uncountable models}
\bigskip
\centerline{
{\bigrm
Saharon Shelah, Heikki Tuuri and Jouko V\"a\"an\"anen
}
\footnote {*}
{
The first author would like to thank the United States--Israel
Binational Science Foundation for support of this research
(Publication \# 377). The second and third author were
supported by Academy of Finland grant 1011040.
}
}

\bigskip
\bigskip

\centerline{\bf Abstract}
\abstract
Let $\s(\A)$ denote the number of automorphisms of a model
$\A$ of power $\o_1$. We derive a necessary and
sufficient condition in terms of trees for the
existence of an $\A$ with $\o_1 < \s(\A) < 2^{\o_1}$.
We  study the sufficiency
of some conditions for $\s(\A) = 2^{\o_1}$.
These conditions are analogous to
conditions studied by D.Kueker in connection with countable
models.

\endabstract

\bigskip
\bigskip
\noindent
The starting point of this paper was an attempt to generalize
some results of D.Kueker [8] to models of power $\o_1$.
For example, Kueker shows that for countable $\A$
the number $\s(\A)$ of automorphisms of $\A$ is either
$\le \o$ or $2^\o$.  In Corollary 13 we
prove the analogue of this result under the set-theoretical
assumption $I(\o)$: if $I(\o)$ holds and
the cardinality of
$\A $  is $\o_1$, then $\s(\A) \le \o_1$ or
$\s(\A) = 2^{\o_1}$.
In Theorem 16 we show that the
consistency strength of this statement + $2^{\o_1} > \o_2$
is that of an inaccessible cardinal.
We use $\dv \A \dv$ to denote the universe
of a model $\A$ and $\ve \A \ve $ to denote the cardinality
of $\dv \A \dv$.
Kueker proves also that if $\ve \A \ve \le \o$,
$\ve \B \ve > \o$ and $\A \equiv \B$ (in $L_{\infty \o}$), then
$\s(\A) = 2^\o$. Theorem 1 below generalizes this
to power $\o_1$.
If $\A$ and $\B$ are countable, $\A \ne \B$
and $\A \prec \B$ (in $L_{\infty\o}$), then we know that
$\s(\A) = 2^\o$. Theorem 7 shows that the natural analogue
of this result fails for models of power $\o_1$.
Theorem 14 links the existence of a model
$\A$ such that $\ve \A \ve = \o_1$, $\o_1 < \s(\A) < 2^{\o_1}$,
to the existence of
a tree $T$ which is of power $\o_1$, of height $\o_1$
and has $\s(\A)$ uncountable branches.

We use $\A \equiv_{\o_1} \B$ to denote
that $\eks$ has a winning strategy in
the Ehrenfeucht-Fra\"\i ss\'e game $G(\A,\B)$
of length $\o_1$ between $\A$ and $\B$.
During this game two players
$\eks$ and $\kai$ extend a countable partial
isomorphism $\p$ between $\A$ and $\B$.
At the start of the game $\p$ is empty.
Player $\kai$ begins the game by choosing
an element $a$ in either $\A$ or $\B$.
Then $\eks$ has to pick an element $b$
in either $\A$ or $\B$ so that $a$
and $b$ are in different models.
Suppose that $a \in \A$.
If the relation $\p \cup \{(a,b)\}$
is not a partial isomorphism,
then $\eks$ loses immediately,
else the game continues in the
same manner and the new value of $\p$ is
the mapping $\p \cup \{ (a,b)\}$.
The case $a \in \B$ is treated similarly,
but we consider the relation $\p \cup \{(b,a)\}$.
The length of our game is $\o_1$ moves.
Player $\eks$ wins, if he can
move $\o_1$ times without losing.
The only difference
between this game and the ordinary
game characterizing partial isomorphism
is its length.
M.Karttunen and T.Hyttinen  have proved
([3,4,7])
that
$\A \equiv_{\o_1} \B$ is equivalent
to elementary equivalence relative to
the infinitely deep language $M_{\infty\o_1}$.
It may also be observed that
$\A \equiv_{\o_1} \B$ is equivalent
to isomorphism
in a forcing extension,
where the set of forcing
conditions is countably closed  [9].
For the definition of $M_{\infty\o_1}$
and other information of $\equiv_{\o_1}$
the reader is referred to
[3,4,7,9,10,11].
Our treatment is selfcontained, however.
The definition of
the language
$M_{\infty\o_1}$ is not needed
in this paper.

One of the basic consequences of
$\A \equiv_{\o_1} \B$ is that if
$\A$ and $\B$ both have power
$\o_1$, then $\A \cong \B$ [7]. The proof of this
is similar to the proof
of the corresponding result for countable models.

We note in passing that there is a canonical
infinitary game sentence $\f_{\A}$
(see [3], [4] or [7]), a kind of
generalized Scott sentence, with
the property that $\B \tar \f_\A$ iff $\A \equiv_{\o_1} \B$
for any $\B$. So, if $\A \equiv_{\o_1} \B$
happens to imply that $\B$ has power $\le \o_1$, then
$\f_\A$ characterizes $\A$ up to isomorphism.

The authors are indebted to Wilfrid Hodges for his
help in the early stages of this work and to
Alistair Lachlan and Alan Mekler for suggesting improvements.

\lause {Theorem 1}. If a model of power
$\o_1$ is $\equiv_{\o_1}$-equivalent
to a model of power $> \o_1$, then
it has $2^{\o_1}$ automorphisms.

\valit
For the proof of this theorem we define
the following game $G(\A)$ where $\A$ is a model of power
$\o_1$:
There are $\o_1$ moves and two players
$\eks$ and $\kai$. During the game
a countable partial isomorphism $\p$ is extended.
At each move $\kai$ first
plays a point, to which
$\eks$ then tries to extend $\p$.
$\kai$ can tell whether the point is to be
on the image side or in the domain side.
Moreover, $\eks$ has to come up with
two contradictory extensions of $\p$,
from which $\kai$ chooses the one the game goes on
with. $\eks$ wins,
if he can play all $\o_1$ moves.

A model $\A$ is called {\it perfect},
if $\eks$ has a winning strategy in $G(\A)$.

\lause {Proposition 2}. If $\A \equiv_{\o_1} \B$
for some $\B$ of power $> \o_1$,
then $\A$ is perfect.

\tod Proof. Let $S$ be a winning strategy of
$\eks$ in the Ehrenfeucht-Fra\"\i ss\'e-game. An
$S$-{\it mapping} is a partial
isomorphism between $\A$ and $\B$ arising
from $S$. We describe a
winning strategy of $\eks$ in
$G(\A)$. During the game $\eks$
constructs $S$-mappings $\s: \A \to \B$
and $\r: \B \to \A$ simultaneously
with the required $\p$. The idea
is to keep  $\p = \r \circ \s$.

\midinsert
\vskip 5cm
\centerline{Figure 1.}
\endinsert

Suppose  now $\kai$ plays
$x$ and asks $\eks$ to
extend the domain of $\p$
to $x$. If $x \notin \dom(\s)$
($=\dom(\p)$), $\eks$ uses
$S$ to extend $\s$ to $x$.
Likewise, if $\s(x) \notin \dom(\r)$,
$\eks$ uses $S$ to extend $\r$ so that
$\s(x) \in \dom(\r)$.
Let $\p(x)=\r(\s(x))$.
This completes the first part
of the move of $\eks$.

For the second part, $\eks $ has
to come up with $\p'$ and $\p''$,
which are contradictory extensions of $\p$.
For any $b \in \B$ $S$ gives some $ s(b) \in \A$.
If $b \notin \ran (\s)$, then  $s(b) \notin \dom (\p)$.
As $\vert \B \backslash \ran(\s) \vert>\vert
\A\vert$, there are $b \ne b' \in \B \backslash \ran(\s)$
with $s(b)=s(b')$.
We extend $\r$ using $S$ first to get an element
$a$ so that $\r(b)=a$ and after that we extend
$\r$ further to get $\r(b')=a'$. Now,
$a \ne a'$, since $b\ne b'$ (Figure 1).
Now we can define $\p'$ and $\p''$.
In the first case we extend
$\s$ so that $\s(s(b))=b$
and we let
$\p'=\r \circ \s$.
(Note here, that we do not extend $\s$ to $b'$.
It is not necessary to keep $\ran(\s) = \dom(\r)$.)
In the second case we extend
$\s$ so that $\s(s(b))=b'$ and
we define
$\p''=\r \circ \s$.
Because $\p'(s(b)) \ne \p''(s(b))$, the two
extensions are contradictory.
$\eop$

\lause {Proposition 3}. If $\A$ is perfect,
then $\s(\A)=2^{\o_1}$.

\tod Proof. Suppose $S$
is a winning strategy of $\eks$ in $G(\A)$. Let us
consider all games in which $\kai$
enumerates all of $\A$. Each such
play determines an automorphism of
$\A$. Since $\kai$ has a chance of splitting
the game at each move,
there are $2^{\o_1}$ different automorphisms.
$\eop$

\valit This ends the proof of Theorem 1. $\eop$

\valit
Now we define a game that characterizes the elementary
submodel relation for the language
$M_{\infty \o_1}$.
Suppose $\A \subseteq \B$. We describe the
game $G_{\preceq}(\A,\B)$.
The game resembles very much the ordinary
Ehrenfeucht-Fra\"\i ss\'e-game between $\A$ and $\B$. The difference is
that at the start of the game $\kai$ can pick
a countable set $C$ of elements of $\A$
and set as the initial partial isomorphism
$\p=\{(a,a) \ \vert \ a \in C \}$.
Then $\kai$ and $\eks$ continue the game
like the usual Ehrenfeucht-Fra\"\i ss\'e-game extending $\p$.

We write $\A \preceq_{\o_1} \B$, if $\eks$ has
a winning strategy in the game
$G_{\preceq}(\A,\B)$.
If $\A \preceq_{\o_1} \B$ and $\A \ne \B$, then we write
$\A \prec_{\o_1} \B$.
It can be proved that
the relation $\A \preceq_{\o_1} \B$ holds if and only if
$\A$ is an elementary submodel of $\B$
relative to the language $M_{\infty \o_1}$.
In this definition the formulas of $M_{\infty \o_1}$ may contain
only a countable number of free variables.
The proof is very similar to the proof
of the fact that $\A \equiv_{\o_1} \B$ is equivalent to
elementary equivalence of $\A$ and $\B$ ([7], [3], [4]).

We describe the game $G_{\le}(\A,\B)$, which
is more difficult for $\eks$ to win than
$G_{\preceq}(\A,\B)$.
The length of the game is $\o_1$ and it
resembles the Ehrenfeucht-\Fraisse game.
During it $\eks$ must extend a countable partial
isomorphism $\pi : \A \to \B$ and
at each move the rules are the
following:
\itskip
\ti if $a \in \A$, $a \notin \dom(\pi)$ and $a \notin \ran (\pi)$,
then $\kai$ can move $a \in \A$ and demand $\eks$ to extend
$\pi$ to $\pi \cup \{ (a,a) \}$;
\tii if $a \in \A$ ($a \in \B$) then $\kai$
can move $a \in \A$  ($a \in \B$) and demand $\eks$ to
extend $\pi$ so that $a \in \dom(\pi)$ ($a \in \ran(\pi)$).
\itskip

We write $\A \le_{\o_1} \B$, if
$\A \subseteq \B$ and
$\eks$ has
a winning strategy in the game $G_\le(\A,\B)$.
If $\A \le_{\o_1} \B$ and
$\A \ne \B$, then we write $\A <_{\o_1} \B$.

Our aim is next to prove that if $\A <_{\o_1} \B$
for some $\B$, then there are $2^{\o_1}$
automorphisms of $\A$.

\lause {Lemma 4}.
Let $(\A_\a)_{\a < \d}$ ($\d$ limit) be
uncountable models such that:
\itskip
\ti $\A_\a \subseteq \A_\b$ if $\a < \b$;
\tii $\A_\g = \bigcup_{\a < \g} \A_\a$ if $\g$ is
a limit;
\tiii $\A_\a \le_{\o_1} \A_{\a+1}$ if
$\a <\d$.
\itskip
\noindent
Let $\A_\d= \bigcup_{\a < \d} \A_\a$.
Then $\A_0 \preceq_{\o_1} \A_\d$.
(The arity of relations and functions must be finite.)

\tod Proof.
For simplicity of notation,
we assume that in the games $G_\preceq(\A,\B)$ and
$G_\le(\A,\B)$
at each round $\a$,
$\eks$ extends the partial isomorphism $\pi$ by
just a single ordered pair $(a_\a,b_\a)$,
where $a_\a \in \A$ and $b_\a \in \B$.

For each $\a < \d$, let $\s_\a$ be $\eks$'s fixed
winning strategy in $G_\le(\A_\a, \A_{\a + 1})$.

We describe a winning strategy for $\eks$ in
$G_\preceq(\A_0, \A_\d)$. We modify the game
$G_\preceq(\A_0, \A_\d)$ so that $\kai$ and $\eks$ only move at
infinite limit ordinal rounds, which is clearly equivalent to the
original game.
At each round $\g < \o_1$, $\eks$ also constructs a sequence $s_\g$ of length
$\d+1$, such that $s_\g(\a) \in \A_\a$ for all $\a \le \d$.
At limit rounds $\g$, $\eks$ first constructs $s_\g$ and then
extends the partial isomorphism $\pi$
in the game $G_\preceq(\A_0, \A_\d)$
by $(a, b)$, where $a = s_\g(0)$ and
$b = s_\g(\d)$.

Before round $\g \ge \o$, we assume that the following
conditions are true:
\itskip

\item{(1)} For all $\a < \d$, the sequence
$((s_\e(\a), s_\e(\a+1)))_{\e < \g}$ is a play
in $G_{\le}(\A_\a, \A_{\a + 1})$
according to $\eks$'s
winning strategy $\s_\a$.

\item{(2)} For all $\e < \g$, $s_\e$ is continuous, that is,
if $\x$ is a limit ordinal, and $s_\e(\x)=a$, there
is $\z < \x$, such that for all $\z < \a  \le \x$,
$s_\e(\a)=a$.

\item{(3)} Suppose $a$ is in the range of some sequence
$s_\e$, $\e < \g$, and $\a$ is the least ordinal such that
$a \in A_\a$. Then there is an ordinal $\b$ such that
$[\a, \b] = \{ \x \mid$ for some $\e < \g$, $s_\e(\x)=a \}$.
If $\g$ is a successor, then $\b$ is a successor ordinal
or $\d$. If $\g$ is a limit, then $\b = \d$.

\itskip
$\kai$ starts the game $G_\preceq(\A_0,\A_\d)$
by choosing the countable set $C$ of elements of $\A_0$.
$\eks$ chooses as the first sequences $s_n$, $n < \o$, constant
sequences whose values enumerate $C$.
Let us consider round $\g$ in the game, where $\g$ is
an infinite limit. In general there are two cases.

First the case where $\kai$ picks $a \in \A_0$
as his $\g$th move. If there is some $s_\e$ such that
$s_\e(0)=a$, then $\eks$ responds by $s_\e(\d) \in \A_\d$
and defines $s_\g = s_\e$.
Else, by (3), $\eks$ can move $a \in \A_\d$ and
choose the appropriate constant sequence as $s_\g$.
The inductive hypotheses are met and we can let $s_{\g+n} = s_\g$,
for $n < \o$.

Suppose then $\kai$ picks $b \in \A_\d$ as his $\g$th move.
Again, if for some $\e< \g$, $s_\e(\d)=b$, we are done.
Else, let us construct the required sequence
$s_\g$. Let $\a_0$ be the least ordinal such that
$b \in \A_{\a_0}$ and $s_\e(\a_0) \ne b$ for all
$\e < \g$. Note that by hypothesis (3)
and condition (ii) of the lemma, $\a_0 = \b_0 + 1$,
for some $\b_0$ (or $\a_0=0$). We define $s_\g(\b) = b$
for all $\b > \b_0$.
Let $c$ be the response of $\eks$ according to $\s_{\a_0}$
if $\kai$ continues
$G_\le (\A_{\b_0}, \A_{\a_0})$ by moving $b \in \A_{\a_0}$.
Let $s_\g(\b_0)=c$. Then we continue the construction
of $s_\g$ by downward induction.
$\eks$ then moves $s_\g(0) \in \A_0$ in the
game $G_{\preceq}(\A_0, \A_\d)$.
Similarly, by a closing procedure, $\eks$
can construct $s_{\g + n}$,
$n< \o$, so that clause (3) is satisfied at $\g + \o$.
$\eop$

\lause {Proposition 5}. If
$\A$ is of cardinality $\o_1$ and
$\A <_{\o_1} \B$ for some $\B$, then
$\A \equiv_{\o_1} \B$ for some
$\B$ of power $\o_2$,
whence $\A$ is perfect.

\tod Proof. We may assume $\A$ and $\B$
have both power $\o_1$.
Thus, by remarks preceding Theorem 1,
$\A \cong \B$.
We construct a sequence
$(\A_\a)_{\a < \o_2}$ of models so that
each is isomorphic to $\A$,
$\A_\a \subset \A_\b$, if $\a < \b$, and
$\A_\a <_{\o_1} \A_{\a+1}$ for all $\a < \o_2$.
We handle the successor step by identifying
$\A_\a$ with $\A$ via the isomorphism.
Then from $\B$ we get $\A_{\a+1}$.
At limits we take the union of models.
Lemma 4 makes sure that the union is
isomorphic to $\A$, if it is not of power $\o_2$.
$\eop$

\valit
So, if $\A$ fulfills the condition of Proposition 5,
then it has $2^{\o_1}$ automorphisms.
The proof of the following result shows that
$\A \le_{\o_1} \B$ is a much stricter
condition than $\A \preceq_{\o_1} \B$.

\lause {Proposition 6}.
$$
\A \le_{\o_1} \B \Rightarrow \A \preceq_{\o_1} \B
$$
but
$$
\A \preceq_{\o_1} \B \not \Rightarrow \A \le_{\o_1} \B.
$$

\tod Proof. The first claim is trivial.
For the second consider the following models.
There is one equivalence relation $R$ in
the vocabulary. The model $\A$
contains simply $\o_1$ equivalence classes of
size $\o_1$.
The model $\B \supset \A$ contains one additional
equivalence class of size $\o_1$.
Then it is very easy to see that $\eks$ wins $G_\preceq(\A,\B)$.
But $\kai$ can win $G_\le(\A,\B)$ in two moves.
First $\kai$ chooses some $b \in \B$, $b \notin \A$.
Let $\pi$ be $\kai$'s response. Let $a \in \A$,
$\A \tar R(a, \pi\inv(b))$, $a \notin \ran(\pi) \cup \dom(\pi)$.
Then $\kai$ demands $\eks$ to map $a$ identically.
$\eop$

\valit If $\A$ and $\B$ are countable, $\A \ne \B$ and $\A \prec \B$
(relative to $L_{\infty \o}$), then $\s(\A)=2^\o$.
This would suggest the analogous conjecture for uncountable
models: if $\ve \A \ve = \ve \B \ve = \o_1$ and
$\A \prec_{\o_1} \B$, then $\s(\A)=2^{\o_1}$.
But this conjecture is false, as the following counterexample
constructed by S.Shelah shows.

\lause  {Theorem 7}.
Let $\k > \o$ be regular. There are models
$\M_1 \subseteq \M_2$, $\M_1 \ne \M_2$,
$\ve \M_1 \ve = \ve \M_2 \ve = \k$,
such that
\itskip
\item{(i)} for every $A \subset \dv \M_1 \dv$,
$\ve A \ve  < \k$, there is an isomorphism
from $\M_2$ onto $\M_1$ which is the identity on $A$;
\item{(ii)} $\s(\M_1) \le \k$.

\tod Remark. Hence $\M_1 \prec_\k \M_2$ but there is no
$\M_3$ such that $\M_1 \equiv_\k \M_3$ and $\ve \M_3 \ve > \k$,
as then $\s(\M_1) = 2^\k$.

\tod Proof. We first define such $\M_1$ and $\M_2$ with
the vocabulary
$L= \{ R_\d \mid 0 < \d < \k, \d {\rm\ limit} \}$,
where $R_\d$ has $\d$ places and $\ve R_\d^{\M_1} \ve
= \ve R_\d^{\M_2} \ve = \k$. We can then replace these
models (in Proposition 8) by models with a vocabulary
consisting of just one binary relation.

We define $A$, $A_\a$, $f^\a$ and $\g_\a$, $\a < \k$,
such that:
\itskip
\item{(1)} $\o \le \g_\a < \k$ for all $\a < \k$ and
$\langle \g_\a \mid \a < \k \rangle$ is increasing and
continuous;
\item{(2)} $\g_0 = \o$,
if $\a> 0$ is a limit, then $\g_\a = \bigcup_{\b < \a}
\g_\a$, and if $\a = \b + 1$, then $\g_\a = \g_\b + \g_\b$;
\item{(3)} $A_\a = \{ i < \g_\a \mid i {\rm \ even} \}$,
$A = \{ i < \k \mid i {\rm \ even} \}$;
\item{(4)} $f^\a$ is a 1--1 function  from $\k$ onto $A$
mapping $\g_{\a+1}$ onto $A_{\a+1}$;
\item{(5)} $f^\a$ maps the interval $[ \g_\b, \g_{\b+1} )$
onto $[ \g_\b, \g_{\b+1} ) \cap A$ for $\b > \a$;
\item{(6)} $f^\a \raj A_\a$ is the identity function on $A_\a$;
\item{(7)} $f^\a$, $\a < \k$, are defined
using free groups (see the construction of $f^\a$ below).
\itskip
The definition of $\g_\a$ and $A_\a$ is clear from (1)--(3).
We now describe the construction of $f^\a$, $\a < \k$.
If $\b < \k$, let $T^\b_{\rm at}= \{ s^\b_\a \mid \a \le \b \}$
and $T^\b_{\rm nat} = \{ (s^\b_\a)\inv \mid \a \le \b \}$ be
sets of arbitrary symbols. Let $T_\b$ be the set of all such sequences
$\t= \s_1 \ldots \s_n$ that:
\itskip
\item{(T1)} $0 \le n < \o$;
\item{(T2)} $\s_k \in T^\b_{\rm at} \cup T^\b_{\rm nat}$
for all $1 \le k \le n$;
\item{(T3)} if $n > 0$ then $\s_n = s^\b_\b$;
\item{(T4)} $\s_k \in T^\b_{\rm nat}$ $\Rightarrow$
 $ \s_{k+1} \in T^\b_{\rm at}$ for all $1 \le k < n$;
\item{(T5)} $\neg (\eks k , \a) ( \{ \s_k, \s_{k+1} \}
    = \{ s^\b_\a, (s^\b_\a)\inv \} )$.
\itskip
Thus we see that $T_\b$ is a subset of
the normal forms of the free group generated
by $\{ s^\b_\a \mid \a \le \b \}$.
If $\t = \s_1 \ldots \s_n \in T_\b$ and $s^\b_\a \in T^\b_{\rm at}$,
then we define the operation $s^\b_\a \cdot \t$ in
the following way:
\itskip
\item{(a)} if $\s_1 \ne (s^\b_\a)\inv$ or $\t=\eset$,
then
$s^\b_\a \cdot \t = s^\b_\a \s_1 \ldots \s_n$ (i.e. just concatenate);
\item{(b)} if $\s_1 = (s^\b_\a)\inv$, then $s^\b_\a \cdot \t =
 \s_2 \ldots \s_n$.
\itskip
\noindent
It is easy to check that $s^\b_\a \cdot \t \in T_\b$.
Thus $\cdot$ is
defined like the multiplicative operation for the free
group.

\sublemma {Lemma A}.  Let $\t, \t' \in T_\b$ and $\a \le \b$.
If $\t \ne \t'$, then $s^\b_\a \cdot \t \ne s^\b_\a \cdot \t'$.

\tod Proof. Straightforward.
$\eop$ Lemma A.

\medskip

For each $\a < \k $ let
$$
\{ (\t_\x, j_\x) \mid \g_\a \le \x < \g_{\a+1} \}
$$
list the set
$$
P_\a = \{ (\t, j) \mid \t \in T_\a, \t \ne \eset,
j < \g_\a, j \notin A_\a \}
$$
without repetitions in such a way that
\medskip

\centerline{$\x$ is even if and only if
$\s^{\t_\x}_1 \in T^\a_{\rm at}$,}
\medskip

\noindent
where we denote $\t_\x = \s^{\t_\x}_1 \ldots \s^{\t_\x}_{n_{\t_\x}}$.

If $(\t,j) \in P_\a$ for some $\a < \k$, let
$\x(\t, j)$ be the unique $\x$ such that $(\t,j) = (\t_\x , j_\x)$.
Now we define $f^\a$, $\a < \k$ (see Figure 3). For $\e < \k$ let
$$
f^\a(\e) = \cases{
\e              &if $\e < \g_\a$ and $\e \in A_\a$, \cr
\x(s^\a_\a, \e) &if $\e < \g_\a$ and $\e \notin A_\a$, \cr
\x(s^\a_\a \cdot \t, j) &if $\g_\a \le \e < \g_{\a+1}$ and
                          $\e = \x(\t,j)$, \cr
\x(s^\b_\a \cdot \t, j) &if $\g_\b \le \e < \g_{\b+1}$,
$\b > \a$ and $\e = \x(\t,j)$. \cr
}
$$

\midinsert
\vskip 3.5cm
\centerline{Figure 3.}
\endinsert

We have to check that $f^\a$ is well-defined, that is,
$\x(s^\a_\a \cdot \t, j)$ and $\x(s^\b_\a \cdot \t, j)$
must be defined above in appropriate conditions and their
values must be even. We check only $\x(s^\a_\a \cdot \t, j)$,
the other case is similar.
Suppose $\g_\a \le \e < \g_{\a+1}$ and $\e=\x(\t,j)$. Then
$\t \in T_\a$, $\t \ne \eset$.
Let $\t=\s_1 \ldots \s_n$. If $\s_1 \ne (s^\a_\a)\inv$, then
$s^\a_\a \cdot \t = s^\a_\a \s_1 \ldots \s_n \ne \eset$. Thus
$\x(s^\a_\a \cdot \t, j)$ is defined and it is even,
since $s^\a_\a \in T^\a_{\rm at}$.
Suppose $\s_1 = (s^\a_\a)\inv$. Then
$s^\a_\a \cdot \t = \s_2 \ldots \s_n$. Now
$n \ge 2$ by (T3) and $\s_2 \in T^\a_{\rm at}$ by (T4).
Thus $\s_2 \ldots \s_n \ne \eset$ and
$\x(s^\a_\a \cdot \t , j)$ is defined and even.

\sublemma {Lemma B}. Conditions (4), (5) and (6) above are met.

\tod {Proof.} From the definition of $f^\a$ we see
easily that $f^\a$ maps $\g_{\a+1}$ to $A_{\a+1}$ and
$[ \g_\b , \g_{\b+1} )$ to $[ \g_\b, \g_{\b+1}) \cap A$,
if $\b > \a$.
We show first that $f^\a$ is a 1--1 function $\k \to A$.
Suppose $\e_1 \ne \e_2$. We prove $f^\a(\e_1) \ne f^\a(\e_2)$.
There are several cases, of which we treat the two most
interesting. The proof in other cases is similar or trivial.

(a) Suppose $\e_1 < \g_\a$, $\e_1 \notin A_\a$
and $\e_2 \in [\g_\a, \g_{\a+1})$. Let $\e_2=\x(\t,j)$.
Since $\t \ne \eset$, by Lemma A $s^\a_\a \cdot \t \ne s^\a_\a$.
Thus $f^\a(\e_1) = \x(s^\a_\a, \e_1) \ne \x(s^\a_\a \cdot \t, j)
= f^\a(\e_2)$.

(b) Suppose $\e_1,\e_2 \in [\g_\a, \g_{\a+1} )$.
Let $\e_1 = \x(\t_1, j_1)$ and $\e_2 = \x(\t_2, j_2)$.
If $j_1 \ne j_2$, then the claim is clear.
If $j_1 = j_2$, then $\t_1 \ne \t_2$ and by Lemma A
$s^\a_\a \cdot \t_1 \ne s^\a_\a \cdot \t_2$ and again the claim holds.

Next we prove that $f^\a$ is onto.
Let $\d \in A$. We try to find $\e < \k$, for which $\d=f^\a(\e)$.
If $\d \in A_\a$, then we set $\e = \d$.
Suppose then $\d \in [\g_\a, \g_{\a+1} ) \cap A$.
Denote $\d=\x(\t,j)$, where $\t= \s_1 \ldots \s_n$, $\t \ne \eset$.
We know $\s_1 \in T^\a_{\rm at}$, since $\d$ is even.
\itskip
\item{(a)} If $n=1$, then   $\t= s^\a_\a$ by (T3) and we set
$\e = j$.
\item{(b)} If $n > 1$ and $\s_1 = s^\a_\a$, then
we set $\e = \x(\s_2 \ldots \s_n, j)$.
\item{(c)} If $n>1$ and $\s_1 \ne s^\a_\a$, then
$\e= \x((s^\a_\a)\inv \s_1 \ldots \s_n, j)$.
Here $\x$ is defined and (T4) fulfilled because $\s_1 \in T^\a_{\rm at}$.
\itskip
\noindent
Suppose then $\d \in [\g_\b, \g_{\b+1}) \cap A$, $\b >\a$.
\itskip
\item{(a)} If $\s_1 = s^\b_\a$, then $n>1$ by (T3)
and $\e = \x(\s_2 \ldots \s_n, j)$.
\item{(b)} If $\s_1 \ne s^\b_\a$, then
$\e = \x((s^\b_\a)\inv \s_1 \ldots \s_n, j)$.
\itskip
Thus we have proved that $f^\a : \k \to A$ is
1--1 and onto. Now (4), (5) and (6) are clear.
$\eop$ Lemma B.
\medskip
If $\a < \k$, let $\g(\a)$ denote the unique $\b$
for which $\g_\b \le \a < \g_{\b+1}$.
Let $G_1$ be the group of permutations of $A$ generated
by $\{ f^\b (f^\a)\inv \mid \a, \b < \k \}$.
Let $G_2$ be the group of permutations of $\k$
generated by $\{ (f^\b)\inv f^\a \mid \a,\b < \k \}$.

We are ready to define the models. We define  $\M_1$ and $\M_2$
as follows:
\itskip
\item{(i)} $\dv \M_1 \dv = A$;
\item{(ii)} $ \dv \M_2 \dv = \k$;
\item{(iii)} $R^{\M_k}_\a = \{ \langle i_0 i_2 \ldots i_\e \ldots \rangle
_{\e < \a, \e {\rm \ even}} \mid
 \eks g \in G_k (\bigwedge_{\e < \a {\rm \ even}}
 g(i_\e) = \e)  \}$, $k = 1,2$, $0 < \a < \k$, $\a$ limit.

\sublemma {Lemma C}.  $\M_1 \subseteq \M_2$.

\tod {Proof.} Suppose $\langle i_\e \mid \e < \a {\rm \ even}\rangle
\in R^{\M_1}_\a$. Thus there are $k < \o$, $\a_r, \b_r < \k$,
for $1 \le r \le k$ such that (using $(f^\b (f^\a)\inv)\inv = f^\a (f^\b)\inv$)
$$
\bigwedge_{\e < \a {\rm \ even}} f^{\b_1} (f^{\a_1})\inv
f^{\b_2}(f^{\a_2})\inv \ldots f^{\b_k} (f^{\a_k})\inv (i_\e) =\e.
$$
If $\g < \k$ is chosen large enough,
then by (6) $f^\g(i_\e)=i_\e$ and $f^\g(\e)=\e$
for all $\e < \a$, $\e$ even, and thus
$$
\bigwedge_{\e < \a {\rm \ even}}
((f^\g)\inv f^{\b_1}) ((f^{\a_1})\inv f^{\b_2})
\ldots ((f^{\a_{k -1}})\inv f^{\b_k})
((f^{\a_k})\inv f^\g) (i_\e) = \e.
$$
But this means $\langle i_\e \mid \e < \a {\rm \ even} \rangle
\in R^{\M_2}_\a$. The other direction is similar.
$\eop$ Lemma C.

\sublemma {Lemma D}. For each $\a$, $f^\a$ is an
isomorphism from $\M_2$ onto $\M_1$ which is the identity
on $A_\a$. (Hence $G_k$ is a group of automorphisms of $\M_k$.)

\tod {Proof.}  Suppose $\langle i_\e \mid \e < \a {\rm \ even} \rangle
\in R^{\M_2}_\a$.
Then
$$
\bigwedge_{\e < \a {\rm \ even}}
(f^{\b_1})\inv f^{\a_1} \ldots (f^{\b_k})\inv f^{\a_k} (i_\e) = \e.
$$
If $\g$ is chosen large enough,
then
$$
\bigwedge_{\e < \a {\rm \ even}}
f^\g (f^{\b_1})\inv f^{\a_1} \ldots (f^{\b_k})\inv f^{\a_k}
(f^\a)\inv (f^\a(i_\e)) = \e,
$$
which means
$\langle f^\a(i_\e) \mid \e < \a {\rm \ even} \rangle
\in R^{\M_1}_\a$.
The other direction is similar.
$\eop$ Lemma D.

\medskip
Since $\k$ is regular,
Lemma D proves part (i) of the theorem. To show
(ii) it is enough to prove the following lemma, because
$\ve G_1 \ve \le \k$.

\sublemma {Lemma E}. $G_1$ is the group of all automorphisms of $\M_1$.

\tod {Proof.}
Let $g^* \in {\rm AUT}(\M_1)$, $g^* \notin G_1$. Let
$G_1^\d$ be the group generated by $\{ f^\b (f^\a) \inv \mid
\a,\b < \d \}$.
As $\k$ is regular,
by taking successive closures we can find
a limit ordinal $\d < \k$ such that:
\itskip
\item{($\d 1$)} $g^*$ maps $A_\d$ onto $A_\d$;
\item{($\d 2$)} for every $g \in G^\d_1$, $g^* \raj A_\d \ne g \raj A_\d$.
\itskip
\noindent
(In fact the set of such $\d$ is a closed unbounded
subset of $\k$.)

Let $i_\e = g^*(\e)$ for $\e < \a$, $\a =\g_\d$, $\e$ even.
As $g^* \in {\rm AUT} (\M_1)$ and
$\langle \e \mid \e < \a {\rm \ even} \rangle \in R^{\M_1}_\a$,
there is some $g_1 \in G_1$ with
$\bigwedge_{\e < \a {\rm \ even}} g_1 (i_\e) = \e$.
Let $g=g_1\inv \in G_1$. Then
$$
\bigwedge_{\e < \a {\rm \ even}} g(\e)=i_\e.
$$
Thus $g^* \raj A_\d = g \raj A_\d$. By ($\d 2$)
$g \raj A_\d \notin \{ h \raj A_\d \mid h \in G^\d_1 \}$
and by ($\d 1$) $g$ maps $A_\d$ onto itself.
To get a contradiction it is enough to prove:
\medskip
\item{($\Gamma$)}  If $g \in G_1$ and
$g \raj A_\d \notin \{ h \raj A_\d \mid h \in G^\d_1 \}$,
then $g$ does not map $A_\d$ onto itself.

\tod {Proof of ($\Gamma$).} So let
$$
g = f^{\b_k} (f^{\a_k})\inv \ldots f^{\b_1} (f^{\a_1})\inv
$$
be a counterexample with $k$ minimal.
Clearly $\a_i \ne \b_i$ and $\a_{i+1} \ne \b_i$  by the minimality
of $k$.

As $g \notin G^\d_1$, for some $1 \le r \le k$ holds
$\a_r \ge \d$ or  $\b_r \ge \d$. If $\a_r \ge \d$,
then we can consider
$g\inv = f^{\a_1} (f^{\b_1})\inv \ldots f^{\a_k} (f^{\b_k})\inv$,
which is also a counterexample with $k$ minimal.
Thus we may assume without loss of generality that
$\b_r \ge \d$ for some $r$.
Let
$$
\m = \max ( \{ \a_r \mid r \in \{ 1, \ldots, k \},
\a_r < \d \} \cup
\{ \b_r \mid r \in \{ 1, \ldots, k \}, \b_r < \d \} ) + 1.
$$

Let $\x_0 \in A_\d$ be arbitrary. We denote
$$
\eqalign{
\n_1 &= (f^{\a_1})\inv(\x_0), \cr
\x_1 &= f^{\b_1}(\n_1), \cr
&\vdots \cr
\n_k &= (f^{\a_k})\inv(\x_{k-1}), \cr
\x_k &= f^{\b_k}(\n_k).
}
$$
Thus $\x_k = g(\x_0)$. For $i=0, \ldots, k$ let
$$
b_{\le i} = \max \{\m, \b_1, \ldots, \b_i \}.
$$

\sublemma {Lemma F}. Suppose $\x_0 \in A_\d$. Then
$\g(\x_i) < \max \{ b_{\le i} + 1, \d \}$ for $i=0,\ldots,k$.

\tod {Proof.} By induction.
First, $\g(\x_0) < \d$. Suppose
$\g(\x_i) < \max \{ b_{\le i} + 1, \d \}$.
\relax From the definition of $f^\a$ we see
 $\g((f^\a)\inv(\e)) \le \g(\e)$ for all
$\a,\e$. Thus  $\g(\n_{i+1}) \le \g(\x_i)$.
We see also that if
$\g(f^\a(\e)) > \g(\e)$, then $\g(f^\a(\e))=\a$.
Thus $\g(\x_{i+1}) \le \g(\n_{i+1})$ or
$\g(\x_{i+1})=\b_{i+1}$. In both cases
$\g(\x_{i+1}) < \max \{ b_{\le i+1} + 1, \d \}$.
$\eop$ Lemma F.

\sublemma {Lemma G}. For all $1 \le i \le k$ either
$\b_i \le b_{\le i -1}$ or $\a_i \le b_{\le i -1}$.

\tod Proof. Suppose $\b_i > b_{\le i -1}$ and
$\a_i > b_{\le i - 1}$.
Since $b_{\le i - 1} \ge \m$, this implies
$\a_i, \b_i \ge \d$.
Thus $\a_i, \b_i \ge \max \{ b_{\le i - 1} + 1, \d \}$.
Suppose $\x_0 \in A_\d$ is arbitrary.
By Lemma F $\g(\x_{i-1})  < \max \{ b_{\le i - 1} + 1, \d \}$
and by (6) $f^{\b_i} (f^{\a_i})\inv (\x_{i-1}) =\x_{i-1}$.
But now we see
$$
\eqalign{
f^{\b_k}&(f^{\a_k})\inv \ldots f^{\b_1}(f^{\a_1})\inv \raj A_\d \cr
&=
f^{\b_k}(f^{\a_k})\inv \ldots f^{\b_{i+1}}(f^{\a_{i+1}})\inv
f^{\b_{i-1}}(f^{\a_{i -1}})\inv \ldots f^{\b_1}(f^{\a_1})\inv
\raj A_\d,
}
$$
a contradiction with the minimality of $k$.
$\eop$ Lemma G.

\medskip
The following lemma shows that $g$ maps $\x(s^\m_\m, 1)$
outside $A_\d$, which contradicts our assumption
and proves ($\G$).

\sublemma {Lemma H}.  Let $\x_0=\x(s^\m_\m, 1)$.
Then for all $1 \le i \le k$
$\x_i$ is of the form
$\x(s^{b_{\le i}}_{\b_i} \s^i_2 \ldots \s^i_{n_i}, j_i)$,
where
$s^{b_{\le i}}_{\b_i} \s^i_2 \ldots \s^i_{n_i}
 \in T_{b_{\le i}}$ and
$n_i \ge 1$.
Hence $\g(\x_i) = b_{\le i}$.

\tod Proof. Suppose first the claim holds for $\x_i$, $i \ge 1$. We
prove it holds for $\x_{i+1}$.

(a) Suppose $\a_{i+1} > b_{\le i} = \g(\x_i)$.
Then $\n_{i+1} = (f^{\a_{i+1}})\inv (\x_i) = \x_i$. By Lemma
G $\b_{i+1} \le b_{\le  i}$.
Now
$$
\x_{i+1}=f^{\b_{i+1}}(\x_i) = \x(s^{\b_{\le i}}_{\b_{i+1}}
s^{\b_{\le i}}_{\b_i}\s^i_2 \ldots \s^i_{n_i} , j_i).
$$
Hence the claim holds for $i+1$.

(b) Suppose $\a_{i+1} \le b_{\le i} = \g(\x_i)$.
Then
$$
\n_{i+1}=(f^{\a_{i+1}})\inv(\x_i) = \x((s^{\b_{\le i}}_{\a_{i+1}})\inv
s^{\b_{\le i}}_{\b_i}\s^i_2 \ldots \s^i_{n_i} , j_i),
$$
where $\a_{i+1} \ne \b_i$ by the minimality of $k$.
Note that $\n_{i+1}$ is odd.
If $\b_{i+1} > b_{\le i}$, then
$$
\x_{i+1} = f^{\b_{i+1}}(\n_{i+1}) = \x(s^{\b_{i+1}}_{\b_{i+1}},\n_{i+1})
$$
and the claim holds for $i+1$.
If $\b_{i+1} \le b_{\le i}$, then
$$
\x_{i+1} = f^{\b_{i+1}}(\n_{i+1}) = \x(s^{\b_{\le i}}_{\b_{i+1}}
(s^{\b_{\le i}}_{\a_{i+1}})\inv
s^{\b_{\le i}}_{\b_i}\s^i_2 \ldots \s^i_{n_i} , j_i),
$$
where $\b_{i+1} \ne \a_{i+1}$ by the minimality of $k$
and the claim holds.

Next we prove that the claim is true for $i=1$.

(a) Suppose $\a_1 > b_{\le 0} = \m$. Then
$\n_1 = \x_0=\x(s^\m_\m,1)$ and $\b_1 \le b_{\le 0} = \m$.
As above we get $\x_1=\x(s^\m_{\b_1} s^\m_\m, 1)$.

(b) Suppose $\a_1 \le b_{\le 0} = \m$. Then
$\n_1=\x((s^\m_{\a_1})\inv s^\m_\m, 1)$, where
$\a_1 \ne \m$ by the definition of $\m$.
If $\b_1 > \m$, then $\x_1 = \x(s^{\b_1}_{\b_1}, \n_1)$.
If $\b_1 \le \m$, then $\x_2=(s^\m_{\b_1} (s^\m_{\a_1})\inv s^\m_\m,1)$.
$\eop$ Lemma H.
\medskip

Let $\x_0= \x(s^\m_\m,1)$. By Lemma H
$\g(\x_k)=b_{\le k} \ge \d$, since
$\b_i \ge \d$ for some $i$.
Thus $\x_k \notin A_\d$, which proves ($\G$).
This ends the proof of Lemma E and the whole theorem.
$\eop$

\lause {Proposition 8}. We can find models
$\M_1$ and $\M_2$ which satisfy Theorem 7 and have
a vocabulary of one binary relation.

\tod Proof.
Suppose $\M$ is a model of the vocabulary
$\{ R_\d \mid 0 < \d < \k, \d {\rm \ limit} \}$, such that
$\ve \M \ve = \k$, $\ve R^\M_\d \ve \le \k$ and $R_\d$
has $\d$ places.
We define a model $\A = F(\M)$ of one binary relation
$R$.
Let
$$
\dv \A \dv = \dv \M \dv \cup \bigcup_{\d}
\{ ((a_\a)_{\a < \d}, \b ) \mid
\M \tar R_\d (a_0, \ldots, a_{\a < \d}, \ldots ),
\b < \d\}.
$$
The relation $R$ holds in $\A$ exactly
in the following two cases:
\itskip
\ti if $b_1, b_2 \in \dv \A \dv$,
$b_1=((a_\a)_{\a < \d}, \b_1)$ and
$b_2=((a_\a)_{\a < \d}, \b_2)$, where $\b_1 < \b_2$,
then $\A \tar R(b_1, b_2)$;
\tii if $b \in \dv \A \dv$ and $b=((a_\a)_{\a < \d}, \b)$, then
$\A \tar R(a_\b,b )$.
\itskip
\noindent In other words, for each tuple $(a_\a)_{\a < \d}$,
such that $\M \tar R_\d(a_0,  \ldots, a_{\a < \d}, \ldots)$
we add $\d$ new elements to $\dv \A \dv$.
The new $\d$ elements are wellordered by $R$ and for all $\b < \d$
$a_\b$ is in relation $R$ with the $\b$th added element.

Obviously $\ve F(\M) \ve =  \k$. It is a routine task
to check that there is a 1--1 correspondence between
${\rm AUT}(\M)$ and ${\rm AUT}(F(\M))$.
(Note that $\A \tar \neg \eks x R(x,a)$ iff $a \in \dv \M \dv$.)
Thus $\s(\M)= \s(F(\M))$.
It also easy to see that if $\M \prec_\k \M'$, then
$F(\M) \prec_\k F(\M')$.
Let $\M$ and $\M'$ be the models constructed in Theorem 7.
Let $\M_1 = F(\M)$ and $\M_2 = F(\M')$.
$\eop$

\valit
We say that a chain of models $(\A_\a)_{\a < \k}$ is
{\it continuous}, if $\A_\g = \bigcup_{\a < \g} \A_\a$ for $\g$  a limit.
A chain is an
{\it elementary} chain, if $\A_\a \preceq_{\o_1} \A_\b$ for all
$\a < \b$.
If the relation $\preceq_{\o_1}$
were preserved under unions of
continuous chains of models, then we could replace
$<_{\o_1}$ by $\prec_{\o_1}$ in
Proposition 5, as is easy to see.
This raised the question, whether $\preceq_{\o_1}$
is preserved under unions of continuous chains.
Since Theorem 7 shows that $<_{\o_1}$ cannot be replaced by $\prec_{\o_1}$,
it also proves that $\preceq_{\o_1}$ is
not always preserved. Below we present also two other
counterexamples.
They are continuous elementary chains of length $\o$
and $\o_1$.
The problem, whether $\preceq_{\o_1}$ is
preserved under unions of continuous chains of length
$\o_2$ or greater, is open to the authors.

We define the linear order $\n$, which we shall use
in the proofs below. The linear order $\n$ consists of
functions $f : \o \to \o_1$, for which the
set $\{ n \in \o \ \vert \ f(n)\ne 0\}$ is finite.
If $f,g \in \n$, then $f < g$ iff $f(n) < g(n)$,
where $n$ is the least number, where $f$ and
$g$ differ.
By $\n^{< \a}$ we mean the restriction of $\n$
to those functions $f$ for which $f(0)< \a$.
Similarly we define $\n^{\ge \a}$.

Let $\x$ and $\h$ be arbitrary linear orders.
By $\x \times \h$ we mean a linear order where we
have a copy of $\x$ for every $x \in \h$.
The order between the copies
is determined
by $\h$. By $\h + \x$ we mean
a linear order, where $\x$ is on top of $\h$.
If $\a$ is an ordinal, then $\a^*$ denotes
$\a$ in a reversed order.

We first prove a lemma about $\n$.

\lause {Lemma 9}.

\ti $\n^{\ge \a} \cong \n$ for all
$\a$,
\tii $\n \times n \cong \n$ for
all $n \in \o$,
\tiii $\n \times \a^* \cong \n$ for
all $\a < \o_1$.

\tod Proof.
(i) Let $f \in \n^{\ge \a}$. Simply map
$f$ to $g \in \n$, where $g(0)=f(0)-\a$ and
$g(n)=f(n)$, if $n \ne 0$.

(ii) We prove the claim by
induction on $n$. Suppose $\n \times n \cong \n$.
Clearly $\n^{< 1} \cong \n$,
thus $\n \times n \cong \n^{< 1}$.
By (i)
$\n \cong \n^{\ge 1}$.
So $\n \times (n+1) \cong \n^{< 1} + \n^{\ge 1} \cong \n$.

(iii)
We prove this by induction on $\a$.
The successor step is easy, because
$\n + \n \cong \n$.
Suppose then that $\a$ is a limit ordinal.
Let $(\a_n)_{n < \o}$ be an increasing
sequence cofinal in $\a$.
Then $\a= \sum_{n < \o} \a_{n+1} - \a_n$.
All the differences in the sum are $< \a$,
so we can use our induction assumption
and we get $\n \times \a^* \cong \n \times \o^*$.
Thus the limit case is reduced to showing that
$\n \times \o^* \cong \n$.
We describe the isomorphism.
First we map the topmost copy of
$\n$ in $\o^* \times \n$ to
$\{ f \in \n \ \vert \ f(0) > 0 \}$.
This mapping goes as in (i).
Then we map the next copy of $\n$
to $\{ f \in \n \ \vert \ f(0)=0, f(1)  > 0 \}$,
and continuing
this way we get an isomorphism.
$\eop$

\lause {Proposition 10}.
There exists an elementary chain
$(\A_n)_{n < \o}$ of models of cardinality
$\o_1$ such that
$$
\A_n \not \preceq_{\o_1} \bigcup_{n < \o} \A_n
$$
for all $n$.

\tod Proof.
We let $\A_n= \n \times n$.
Then  the union of the chain is
$\A= \n \times \o$.
We can choose an increasing sequence
of points in $\A$ so that
the length of the sequence is $\o$
and the sequence has no upper bound in $\A$.
It is not possible to find such a
sequence in any  $\A_n$. Thus it is
clear that no  $\A_n$ is an elementary submodel
of $\A$.

It remains to prove that our chain is really an
elementary chain.
We start to play the game $G_\preceq (\A_n, \A_m)$,
$m > n$.
First $\kai$ chooses a countable set $C$
in $\A_n$, which is mapped identically to $\A_m$.
Some of the points of $C$ are in the topmost
copy of $\n$ in $\A_n$. Let $\a < \o_1$ be so
big that none of these points $f$ has $f(0)\ge \a$.
We form an isomorphism between $\A_n$ and
$\A_m$ so that it maps the points in $C$ identically.
We map the part $\n \times (n-1) + \n^{< \a}$
in $\A_n$ identically to $\A_m$.
The remaining part of $\A_n$ is
$\n^{\ge \a}$ and thus isomorphic to $\n$.
The remaining part of
$\A_m$ is isomorphic to $\n + \n \times (m-n)$
and thus isomorphic to $\n$. So we get an isomorphism
between the remaining parts.
Now $\eks$ can win the game simply by playing
according to our isomorphism.
$\eop$

\lause {Proposition 11}.
There exists an elementary chain
$(\A_\a)_{\a < \o_1}$ of models of cardinality $\o_1$
such that
$$
\A_\a \not \preceq_{\o_1} \bigcup_{\a < \o_1} \A_\a
$$
for all $\a$. In this chain $\A_\g= \bigcup_{\a < \g} \A_\a$,
if $\g$ is a limit ordinal.

\tod Proof.
We let $\A_\a = \n + \n \times \a^*$.
Then there is a descending $\o_1$-sequence
in $\A= \bigcup_{\a < \o_1} \A_\a$, but no
descending $\o_1$-sequence in any
$\A_\a$. This shows that
$\A_\a \not \preceq_{\o_1} \A$.

We have to prove that our chain is elementary.
We start to play the game $G_\preceq ( \A_\a, \A_\b)$,
where $\a < \b$.
First $\kai$ chooses a countable set $C$ of points
in $\A_\a$.
Let $\d < \o_1$ be so big that for no $f \in C$
$f(0) \ge \d$.
We form an isomorphism between our models so
that it maps the points in $C$ identically.
First we map the part $\n \times \a^*$ in
$\A_\a$ identically to $\A_\b$.
We map the part $\n^{< \d}$ in the bottom copy
of $\n$ in $\A_\a$ again identically to $\A_\b$.
Now it remains to map
$\n^{\ge \d}$ to $\n^{\ge \d} + \n \times \g^*$,
where $\g=\b - \a$.
But, according to Lemma 7 (i) and (ii), these both are
isomorphic to $\n$, so we get the isomorphism
between $\A_\a$ and $\A_\b$.
Then $\eks$ wins the game by playing according to
this isomorphism.
$\eop$

\valit We shall now consider a totally different
kind of condition which also guarantees
perfectness.
Let $I(\o)$ denote the assumption (taken from [2])
that
\highskip
\highlight

``there is an ideal $I$ on $\o_2$ which is
$\o_2$-complete, normal,
contains all singletons $\{ \a \}$, $\a < \o_2$, and
$$
I^+ =
\{ X \subseteq \o_2 \ \vert \ X \notin I \}
$$
has a
dense subset $K$ such that every
descending chain of length $< \o_1$ of elements
of $K$ has a lower bound in $K$. ''

\endhighlight

\tod Remark. $I(\o)$ implies that $I$
is precipitous and hence that $\o_2$ is measurable
in an inner model. On the other hand, if
a measurable cardinal is Levy-collapsed to
$\o_2$, $I(\o)$ becomes true [1].

We prove that $I(\o)$ implies CH. Suppose
$2^\o \ge \o_2$. Let $T$ be a full binary
tree of height $\o + 1$. Let $A \subseteq \{ t \in T \mid
{\rm height}(t) = \o \}$, $\ve A \ve = \o_2$.
Let $I$ be the ideal on $A$ given by $I(\o)$.
Now it is very easy to construct
$t_0 < \cdots < t_n < \cdots$ and
$X_0 \supseteq  \cdots \supseteq X_n \supseteq \cdots$,
$n < \o$, such that ${\rm height}(t_n) = n$,
$X_n \in K$, and for all $a \in X_n$ holds $a > t_n$. Now
$\bigcap_{n < \o} X_n$ contains at most one element,
a contradiction.

\lause {Theorem 12}. Assume
$I(\o)$. If a model $\A$ of power
$\o_1$ satisfies $\s(\A) > \o_1$, then
$\A$ is perfect.

\tod Proof. (Inspired by [2].) Let $I$ satisfy $I(\o)$.
We may assume $I$ is an ideal on a set
${\rm AUT}$ of automorphisms of power
$\o_2$.
We describe a winning strategy of $\eks$ in
$G(\A)$. Let $X \subseteq {\rm AUT}$ and
$f \in X$. We say that $f$ is
an $I$-{\it point} of $X$, if for all
countable $\p \subseteq f$, it holds that
$[\p] \cap X \in I^+$, where
$[\p]=$ the set of all extensions of $\p$.

\tod Claim: Every $X\in I^+$ has an $I$-point.
\medskip
\noindent
Otherwise every $f \in X$ has a $\p_f \subseteq f$
with $X \cap [\p_f] \in I$.
Because CH holds, there are only $\o_1$
countable $\p$.
This implies
$X \subseteq \bigcup_{f \in X} X \cap [\p_f] \in I$,
a contradiction.
\medskip
The idea of $\eks$ is to construct a
descending sequence $(X_\a)_{\a < \o_1}$
of elements of $K$. We denote
by $\p_\a$ the countable partial isomorphism
at stage $\a$. The descending sequence is
chosen so that for all
$f \in X_\a$ holds $\p_\a \subset f$.

Suppose the players have played $\a$ moves.
Then $\kai$ demands $\eks$ to extend
$\p_\a$ to a point $x$ and give two contradictory
extensions.
For example, $\kai$ demands $x$ to be
on the domain side.
Because functions $f$ can have only
$\o_1$ different values at $x$ and
$I$ is $\o_2$-closed, we can find
$Y \in I^+$, $Y \subseteq X_\a$, such that all the functions in
$Y$ agree at $x$.
Now let $f$ be an $I$-point of
$Y$ and let $f'$ be an $I$-point of
$Y \backslash \{ f \}$.
Because $f$ and $f'$ are two different
mappings, we can choose countable $\p \subset f$
and $\p' \subset f'$ so that $\p$ and
$\p'$ are contradictory extensions of $\p_\a$
and they are defined at $x$.
Now we can choose $X\in K$ and $X'\in K$, ($X,X' \subseteq Y$),
so that  for all $g \in X$ $\p \subset g$
and for all $g \in X'$ $\p' \subset g$.
The extensions $\p$ and $\p'$ are the demanded
contradictory extensions.
For example,
if $\kai$ picks
$\p$, then we set $X_{\a+1} = X$ and $\p_{\a+1}=\p$.

Limit steps in the game do not cause trouble,
because countable descending chains in $K$
have a lower bound in $K$.
$\eop$

\lause {Corollary 13}. Assume
$I(\o)$.
Then the following condition $(*)$ holds:
\item{$(*)$}
If $\A$ is a model of power $\o_1$, then
the conditions
\itskip
\itemitem{(i)} $\s(\A) > \o_1$,
\itemitem{(ii)} $\s(\A) = 2^{\o_1}$,
\itemitem{(iii)} $\A$ is perfect,
\itskip
\item{} are equivalent.

\tod Remark. T. Jech has proved
[5] it consistent that
$2^{\o}=\o_1$, $2^{\o_1} > \o_2$ and
there is a tree of power
$\o_1$ with $\o_2$ automorphisms.
Hence $(*)$ cannot hold
without some set-theoretical
assumption.
We shall later show that the consistency strength of
$(*)$ is that of an inaccessible cardinal.
Note that $(*)$ implies CH.

\valit
The following result of S.Shelah shows a dependence between trees
and the number of automorphisms of an uncountable model.

\lause {Theorem 14}. Suppose that there exists a tree
$T$ of height $\o_1$ such that:
\itskip
\ti $T$ has $\l$ uncountable branches, where
$\o_1 < \l < 2^{\o_1}$;
\tii each level
in the tree has $\le \o_1$
nodes.
\itskip
\noindent
Then we can build a structure $\M$ of cardinality
$\o_1$ with exactly $\l$ automorphisms.

\tod Proof.
Let $T_\a = \{ t \in T \mid \height(t)=\a \}$ and
$$
G_\a = \{ X \subset T_\a \mid \ve X \ve < \o \}
$$
for each $\a < \o_1$. If $X, Y \in G_\a$, we define
$$
X + Y = (X \backslash Y) \cup (Y \backslash X),
$$
i.e. $X+Y$ is the symmetric difference of $X$ and $Y$.
Clearly, $+$ makes $G_\a$ into an Abelian group. Actually, $G_\a$
is a linear vector space over the field $Z_2 = \{0, 1 \}$,
but below we need only to know that $G_\a$ is Abelian.

Let $G$ be the Abelian group, which consists of all
functions ($\o_1$-sequences) $s: \o_1 \to \bigcup_{\a < \o_1} G_\a$, where
$s(\a) \in G_\a$, and addition is defined coordinatewise:
$(s_1 + s_2)(\a) = s_1(\a) + s_2(\a)$.
If $B=(t_\a)_{\a < \o_1}$ is an
$\o_1$-branch in $T$, then $B$ determines naturally
a sequence $b\in G$, where $b(\a) = \{ t_\a \}$.
Let $G' \subseteq G$ be the Abelian group generated by
all sequences $b$ corresponding to $\o_1$-branches.
(Equivalently, $G'$ is the vector subspace spanned by such sequences.)

Suppose $s \in G'$ is arbitrary.
Then $s = b_1 + \cdots + b_n$ for some
$\o_1$-branches $b_1, \ldots, b_n$.
Clearly, if $t \in T_\a$, then $t \in s(\a)$ iff an odd number
of branches $b_1 , \ldots, b_n$ passes through $t$.
\relax From this we see that if $\a < \b$ and $t \in T_\a$, then
\itskip
\item{$(*)$} $t \in s(\a)$ iff $t$ has an odd number of
successors in $s(\b)$.
\itskip

Let $\M'$ be a model of vocabulary $\{ R_s \mid s \in G' \}$
such that
\itskip
\ti $\dv \M' \dv = \{ s \mid s \in G' \}$;
\tii $\M' \tar R_s(s_1,s_2)$ iff $s_2 = s_1 + s$.
\itskip
\noindent
The model $\M'$ is like an affine space, where the set of points
is $\dv \M' \dv$ and the space of differences
$G'$ is kept rigid.
Obviously, $\ve \M' \ve = \l$
and $\AUT(\M')$ consists
of all mappings $\pi'_s$, $s \in \dv \M' \dv$, where
$\pi'_s(x) = x + s$. Thus $\M'$ has exactly $\l$ automorphisms.

Let $\M$ be a model such that:
\itskip
\ti $\dv \M \dv = \{ s \raj \a \mid s \in \dv \M' \dv, \a < \o_1 \}$;
\tii the vocabulary of $\M'$ is $\{ F \} \cup \{R_s \mid s \in \dv \M \dv \}$;
\tiii $\M \tar R_s (s_1, s_2)$ iff the domains of $s,s_1,s_2$ are
equal and $s_2 = s_1 + s$ (where the sum is defined coordinatewise);
\tiv $\M \tar F(s_1, s_2)$ iff $s_1$ is an initial segment
 of $s_2$.
 \itskip
\noindent
Since $\ve T \ve = \o_1$, there are only $\o_1$ countable
initial segments of $\o_1$-branches, and
$\ve \M \ve = \o_1$. We show that there is a 1--1
correspondence between $\AUT(\M')$ and $\AUT(\M)$.
Let $s
\in \dv \M' \dv$ be arbitrary. Then $\pi'_s \in \AUT(\M')$.
We define from $\pi'_s$  an automorphism $\pi_s$ of
$\M$: if $r \in \dv \M \dv$ and $\dom(r)=\a$, then
$\pi_s (r) = r + s \raj \a$. Obviously,  if $s \ne s'$, then
$\pi_s \ne \pi_{s'}$.

Suppose then $\pi$ is an automorphism of $\M$.
We denote by $s_\eset^\b$
a function,
such that $\dom(s^\b_\eset) = \b$ and $s^\b_\eset(\a) = \eset$
for all $\a < \b$.
We define $s \in G$ in the following way:
$s \raj \b = \pi(s^\b_\eset)$ for all $\b < \o_1$.
We show that $s \in \dv \M' \dv$.
By $(*)$
$\ve s(\a) \ve \ge \ve s(\b) \ve$ if $\a \ge \b$.
Since $\ve s(\a) \ve $ is finite for all $\a$,
there must be $n$ and $\b$ such that
$\ve s(\a) \ve = n$ for all $\a \ge \b$.
Thus from $(*)$ we  see that from $\b$ up
$s$ determines some $\o_1$-branches $b_1, \ldots, b_n$,
such that $s \raj (\o_1 \backslash \b) = b \raj (\o_1 \backslash \b)$,
where $b = b_1 + \cdots + b_n$.
It remains to show that $ s \raj (\b+1) = b \raj (\b + 1)$.
We know $s \raj (\b+1) = \pi(s^{\b+1}_\eset) = s' \raj (\b+1)$
for some $s' \in \dv \M' \dv$. Since $s'(\b) = b(\b)$, $(*)$
implies that $s' \raj (\b + 1) = b \raj (\b +1)$,
and thus $ s = b \in \dv \M' \dv$.

Now it is very easy to show that $\pi = \pi_s$.
Thus there is a 1--1 correspondence and
$\M$ has exactly $\l$ automorphisms.
$\eop$

\tod Remark. If the tree $T$ above is a Kurepa
tree, then the resulting model
$\M$ is clearly not perfect.

\valit
We can modify the preceding proof to get
a suitable model with a finite vocabulary.
We add to the model $\M$ the set
$\{a_s \mid s \in \dv \M \dv\}$ of
new elements and
wellorder them with a new relation $<$.
Then
we can use these new elements to code the relations
$R_s$
into a single relation
and we get a finite vocabulary. This modification does not
affect the number of automorphisms.

Theorem 14 is of use only, if the conditions
in it are consistent with ZFC.
We show that this is indeed the case.

A tree $T$ is a {\it Kurepa tree} if:
\itskip
\ti ${\rm height} (T) = \o_1$;
\tii each level of $T$ is at most countable;
\tiii $T$ has at least $\o_2$ uncountable branches.
\itskip
It is well-known (see e.g. [6]) that
Kurepa trees exist in the constructible universe.
Let $\M$ be a countable standard model of
ZFC + $V=L$.
Let $T$ be a Kurepa tree in $\M$.
Let $\l$ be the number of uncountable branches in $T$.
Now we use forcing
to get a model where $2^{\o_1}>\l$.
We utilize  Lemma 19.7 of [6].
In $\M$ the equation $2^{<\o_1}=\o_1$ holds.
Let $\k>\l$ be such that $\k^{\o_1}=\k$.
Let ${P}$ be the set of all functions $p$
such that:
\itskip
\ti ${\rm dom} (p) \subseteq \k \times \o_1$
and $\vert {\rm dom} (p)\vert < \o_1$,
\tii ${\rm ran} (p) \subseteq \{0,1\}$,
\itskip
\noindent
and let $p$ be stronger than $q$ iff $p \supset q$.
The generic extension $\M[G]$ has the same cardinals
as $\M$ and $\M[G] \tar 2^{\o_1}=\k$.
${ P}$ is a countably closed notion of forcing.
Hence Lemma 24.5 of [6] says that the Kurepa tree
$T$ contains in $\M[G]$ just those branches
that are in the ground model. Thus there are
exactly $\l$ uncountable branches in $T$ also in the
extended model $\M[G]$.
CH is true in $L$, therefore $\M[G] \tar 2^\o = \o_1$
by the countable closure of ${P}$.
We have obtained a model $\M[G]$ of
ZFC + CH with a tree $T$, which has the properties
(i)--(ii) of Theorem 14.

\relax From Theorem 14 and the above remarks
we obtain
a new proof of
Jech's result [5]:
\highskip
\highlight

If ZF is consistent,
then ZFC + $2^\o=\o_1$ +
``there exists a model of cardinality $\o_1$
with $\l$ automorphisms, $\o_1 < \l < 2^{\o_1}$ ''
is consistent.

\endhighlight

\valit If we assume CH, we can prove the other
direction in Theorem 14.

\lause {Proposition 15}.
Assume CH. Suppose that we have a model
$\M$ of cardinality $\o_1$ and $\M$ has $\l$ automorphisms,
$\o_1 < \l < 2^{\o_1}$.
Then there exists a tree $T$ of height $\o_1$
such that the conditions (i)--(ii) in Theorem
14 hold.

\tod Proof.
To avoid some complications, we assume that $\M$ has a
relational vocabulary. If not, we can transform
the vocabulary to relational and that does not
affect the number of automorphisms.
The tree $T$ will consist of partial automorphisms
of $\M$.
Let $(a_\a)_{\a < \o_1}$ enumerate $\M$.
Let $\M_\a=\M \vert \grave {} \ \{a_\b \ \vert \  \b < \a \}$.
We let $T= \{ f \ \vert \ f$ is an automorphism
of some $\M_\a \}$. If $f,g \in T$, then $f \le g$ iff
$g$ extends $f$.

Suppose $f$ is an automorphism of $\M$.
Let $\a < \o_1$ be arbitrary. It may be
that the restriction of $f$ to $\M_\a$ is not a bijection
from $\M_\a$ to $\M_\a$, but by taking
successively closures we find
$\b> \a$, for which $f$ gives an automorphism
of $\M_\b$. Thus $f$ determines an uncountable
branch in $T$.

For the other direction, if we have an uncountable
branch in $T$, it is clear that it determines
an automorphism of $\M$.
Thus $T$ has $\l$ uncountable branches.

The tree $T$ may contain at most
$\o_1 \times \o^\o$ nodes. Since we assumed
CH, this is equal to $\o_1$.
So, each level of $T$ contains $\le \o_1$ nodes.
$\eop$

\lause {Theorem 16}. CH + $(*)$ is equiconsistent
with the existence of an inaccessible cardinal.
Also CH + $2^{\o_1} > \o_2$ + ``for all $\A$ of power
$\o_1$, $\s(\A) > \o_1$ implies $\s(\A) = 2^{\o_1}$''
is equiconsistent with the existence of an
inaccessible cardinal.

\tod Proof. Let $\l$ be a strongly inaccessible
cardinal and $\m \ge \l$ so that
$\m = \m^{\aleph_1}$.
Let ${P} = { Q} \times R$,
where ${Q}$ is the
Levy collapse of $\l$ to $\aleph_2$ (see [6], p. 191)
and ${ R}$ is the set of Cohen conditions for adding
$\m$ subsets to $\aleph_1$.
We show that $V^{\ P} \tar (*)$.
Suppose $p \models \s(\A) > \o_1$.
We may assume, without loss of generality, that
$\A \in V$.
Hence there is a ${\ P}$-name $\tilde f$ and
$p \in {\ P}$ so that
$p \models$ ``$\tilde f$ is an automorphism of
$\A$ and $\tilde f \notin V$.''
For any extension $q$ of $p$ let
$$
f^q = \{ (\a,\b) \mid q \models \tilde f(\a) = \b \}.
$$
Now for each extension $q$ of $p$ and for all
countable sets $A,B \subseteq \o_1$ there are
extensions $q^0$ and $q^1$ of $q$ in ${\ P}$ and
an element $a$ of $\o_1$ so that
\itskip
\ti $A \cup \{ a \} \subseteq \dom(f^{q^0}) \cap
                              \dom(f^{q^1})$,
\tii $B \subseteq \ran(f^{q^0}) \cap
                              \ran(f^{q^1})$,
\tiii $f^{q^0}(a) \ne f^{q^1}(a)$.
\itskip
\noindent Using this fact it is easy to see that
$p \models$ ``$\eks$ wins $G(\A)$''.
This ends the proof of one half of the claims.

For the other half of the first claim we assume that
CH + $(*)$ holds. If $\aleph_2$ is not inaccessible
in $L$, then there is a Kurepa tree with
$\ge \aleph_2$ branches, and hence by
the remark after Theorem 14, a non-perfect model
of cardinality $\o_1$ with $> \o_1$ automorphisms.

For the other half of the second claim
we show that under our assumption
$\aleph_2$ has to be inaccessible in $L$.
For this end, suppose $\aleph_2$ is not inaccessible
in $L$. Then there is $A \subseteq \o_1$ so that
$\aleph_2^{L[A]} = \aleph_2$,
$\aleph_1^{L[A]} = \aleph_1$
and GCH holds in $L[A]$ (see, e.g., Jech [6], p.252).
We shall construct a tree with $\aleph_1$ nodes and
exactly $\aleph_2$ branches.
Let $C$ be the set of $\d$
with $\o_1 < \d < \o_2$,
and $L_\d[A] \tar$ ZFC- + ``there is cardinal $\o_1$ and there are
no cardinals $> \o_1$''. Note that $C \in L[A]$.

If $\g < \b$, we denote by $(L_\b[B],\g)$ a model
of vocabulary $(\in, U_1, U_2)$, where
$U_1$ and $U_2$ are unary relations,
the interpretation of $U_1$ is $B$
and the interpretation of $U_2$ is the single element
$\g \in L_\b[B]$.

We form the Skolem hulls in this proof by choosing
as a witness
the element which is the smallest possible
in the canonical
well-ordering of the corresponding model.

{\bf Fact A.}
An easy argument shows that if $\d \in C$ and
$\g  < \d$, then there cannot
be any gaps between
ordinals which are included in the Skolem
hull of $\o_1 \cup \{ \g \}$
(or $\o_1$, as $\g$ is definable in the model)
in $(L_\d[A],\g)$.

Let $\B$ be the class of pairs
$(\a, (L_\b[B],\g)) \in L[A]$,
where
$L_\b[B] \tar $ ZFC- + ``there is cardinal $\o_1$ and
there are no cardinals $> \o_1$'',
$B = A \cap \o_1^{L_\b[B]}$,
$\a < \o_1^{L_\b[B]}$,
$\g < \b$
and $\g > \o_1^{L_\b[B]}$.

We define a partial ordering of these
pairs as follows:
$$
(\a, (L_\b[B], \g)) < (\a', (L_{\b'}[B'], \g'))
$$
if $\a < \a'$, $\b \le \b'$ and $(L_\b[B],\g)$ is the
transitive collapse of the Skolem hull
of $\a \cup \{\g'\}$ in $(L_{\b'}[B'], \g')$.
We define a tree $T$ as follows:
Nodes of the tree are pairs
$(\a, (L_\b[B], \g)) \in \B$ with
$\a < \b < \o_1$. The ordering
of $T$ is the same as that of $\B$.
The cardinality
of $T$ is $\aleph_1$.

If $G=(\a_\x, (L_{\b_\x}[B_\x],
\g_\x))$, $\x < \o_1$,
is an uncountable branch in $T$, then
the direct limit of
$(L_{\b_\x}[B_\x], \g_\x)$, $\x < \o_1$,
is isomorphic to
some $(L_\d[A], \g)$, where $\d \in C$.
If we denote by $H_\a$ the transitive collapse
of the Skolem hull of $\a \cup \{\g\}$, $\a < \o_1$, in
$(L_\d[A], \g)$, then $(\a, H_\a)$,
$\a < \o_1$, is a branch $H$ in $T$.
A straightforward argument shows that $G$ and $H$ coincide.
So the original branch $G$ is in fact in $L[A]$.
Since $T$ has at most $\aleph_2$ uncountable branches
in $L[A]$, it has at most $\aleph_2$
uncountable branches altogether.
On the other hand, by Fact A above,
$T$ clearly has
at least $\aleph_2$ uncountable branches.
We have
shown that $T$ has
$\aleph_1$ nodes and
exactly $\aleph_2$ uncountable branches.
$\eop$

\valit
In this paper we have considered models of cardinality
$\o_1$ and games of length $\o_1$.
When we generalize the model theory
of countable models to uncountable cardinalities, many problems arise.
We chose to concentrate our
attention on $\o_1$, because it offers the simplest example
of an uncountable cardinal, and
even this simple case seems to present enough problems.
Naturally, the results in this paper can be generalized to many
other cardinalities $\k$, i.e. we can consider models
of power $\k$ and games of length $\k$.
The results 1--6 above are valid for any uncountable
cardinal $\k$. Proposition 10 can be generalized for
any regular uncountable cardinal $\k$,
thus we get an elementary chain
of length $\o$, for which $\preceq_\k$ is not
preserved under the union.
\relax From the ideas of Proposition 11 we obtain the
following result: if $\k$ is a regular uncountable cardinal,
$\l$ is a successor cardinal and $\l \le \k$, then
there is an elementary chain of
length $\l$, for which $\preceq_\k$ is not
preserved under the union.
Theorem 14, which shows a dependence
between trees and automorphisms, holds for any
uncountable $\k$.
Proposition 15 has a counterpart
for any regular uncountable $\k$.

\chapter{References}

\item{[1]}
F.Galvin, T.Jech, M.Magidor,
An ideal game, J. Symbolic Logic vol. 43 no. 2
(1978) 284--292.
\item{[2]}
W.Hodges and S.Shelah, Infinite games and
reduced products,
Ann. Math. Logic 20 (1981) 77--108.
\item{[3]}
T.Hyttinen, Games and infinitary languages,
Ann. Acad. Sci. Fenn. Ser. A I Math. Diss. vol 64 (1987).
\item{[4]}
T.Hyttinen, Model theory for infinite quantifier languages,
Fund. Math. 134 (1990) 125--142.
\item{[5]}
T.Jech, Automorphisms of $\o_1$-trees,
Trans. A.M.S. 173 (1972)
57--70.
\item{[6]}
T.Jech, Set theory  (Academic Press, 1978).
\item{[7]}
M.Karttunen, Model theory for infinitely deep languages,
Ann. Acad. Sci. Fenn. Ser. A I Math. Diss.
vol 50 (1984).
\item{[8]}
D.Kueker, Definability, automorphisms and infinitary
languages, in: Barwise, ed., The syntax and semantics
of infinitary languages (Springer, 1968).
\item{[9]}
M.Nadel and J.Stavi, $L_{\infty \l}$-equivalence,
isomorphism and potential isomorphism,
Trans. A.M.S. 236 (1978)
51--74.
\item{[10]}
J.Oikkonen, How to obtain interpolation for $L_{\k^+\k}$,
in: Drake, Truss, eds., Logic Colloquium '86 (North-Holland, 1988).
\item{[11]}
J.V\"a\"an\"anen, Games and trees in infinitary logic:
A survey, in Quantifiers (eds. M.Krynicki, M.Mostowski
and L.Szczerba), to appear.

\bigskip
Saharon Shelah

Department of Mathematics

Hebrew University of Jerusalem

Jerusalem

Israel

\bigskip

Heikki Tuuri

Department of Mathematics

University of Helsinki

Hallituskatu 15

00100 Helsinki

Finland

\bigskip

Jouko V\"a\"an\"anen

Department of Mathematics

University of Helsinki

Hallituskatu 15

00100 Helsinki

Finland

\loppu